# Three–dimensional compact manifolds and the Poincare conjecture.

Alexander A. Ermolitski

*Cathedra of Mathematics, MSHRC,
st. Nezavisimosti 62, Minsk, 220050, Belarus
E-mail: ermolitski@mail.by*

___

**Abstract:** The aim of the work is to prove the following main theorem.

**Theorem**. *Let $M^3$ be a three–dimensional, connected, simple connected, compact, closed, smooth manifold and $S^3$ be the three–dimensional sphere. Then the manifolds $M^3$ and $S^3$ are diffeomorphic.*

___

## 0 Introduction

We can fix some Riemannian metric $g$ on a manifold $M^n$ of dimension $n$ which defines the length of arc of a piecewise smooth curve and the continuous function $\rho(x; y)$ of the distance between two points $x, y \in M^n$. The topology defined by the function of distance (metric) $\rho$ is the same as the topology of the manifold $M^n$, **[5]**.

We should mention that it will suffice to prove that $M^3$ and $S^3$ are homeomorphic since the existence of a homeomorphism between $M^n$ and $S^n$ ($n=\dim M^n$, $n\leq 6$, $n\neq 4$) implies the existence of a diffeomorphism between them. If $n=7$ then there exist such 28 smooth manifolds that every one from them is homeomorphic to $S^7$ but any two from them are not diffeomorphic.

The proof of the main theorem is based on some notions from **[1]**, **[2]** and that will be considered step by step in the following sections. Some results can be useful in the case when $M^3$ is not simply connected or can be generalized for manifolds of dimension $n>3$.

In section 1, using a smooth triangulation and a Riemannian metric we see that every compact, connected, closed manifold $M^n$ of dimension $n$ can be represented as a union of a $n$–dimensional cell $C^n$ and a connected union $K^{n-1}$ of some finite number of $(n-1)$ — simplexes of the triangulation. A sufficiently small closed neighborhood of $K^{n-1}$ we call a *geometric black hole.* In dimension 3 we have $M^3 = C^3 \cup K^2$.

In section 2, we get some technical results permitting to retract 2–dimensional and 1–dimensional simplexes from $K^2$ having boundaries *i. e.* to obtain a decomposition $M^3 = \widetilde{C}^3 \cup \widetilde{K}^2$ where $\widetilde{K}^2$ contains less simplexes than $K^2$ does.

In section 3, we consider the proof of the main theorem consisting of the realization of several algorithms. The number of 2–dimensional simplexes of the complex $K^2$ becomes less every step and finally we have a decomposition $M^3 = C^3 \cup K^1$ where $K^1$ is a connected and simply connected union of some 1–dimensional simplexes *i. e.* $K^1$ is a tree. Using the section 2 we can retract complex $K^1$ to a point $x_0$ therefore a decomposition $M^3 = C^3 \cup \{x_0\}$ is obtained and $M^3$ is homeomorphic to sphere.

## 1. On extension of coordinate neighborhood

**1°.** Let $M^n$ be a connected, compact, closed and smooth manifold of dimension $n$ and $C^n$ be a cell (coordinate neighborhood) on $M^n$. A standard simplex $\Delta^n$ of dimension $n$ is the set of points $x=(x_1, x_2, ..., x_n) \in \mathbf{R}^n$ defined by conditions

$$0 \leq x_i \leq 1,\ i=\overline{1,\ n},\ x_1+x_2+...+x_n \leq 1.$$

We consider the interval of a straight line connected the center of some face of $\Delta^n$ and the vertex which is opposite to this face. It is clear that the center of $\Delta^n$ belongs to the interval. We can decompose $\Delta^n$ as a set of intervals which are parallel to that mentioned above. If the center of $\Delta^n$ is connected by intervals with points of some face of $\Delta^n$ then a subsimplex of $\Delta^n$ is obtained. All the faces of $\Delta^n$ considered, $\Delta^n$ is seen as a set of all such subsimplexes. Let $U(\Delta^n)$ be some open neighborhood of $\Delta^n$ in $\mathbf{R}^n$. A diffeomorphism $\varphi : U(\Delta^n) \to M^n (\delta^n = \varphi(\Delta^n))$ is called a singular $n$–simplex on the manifold $M^n$. *Faces, edges, the center, vertexes* of the simplex $\delta^n$ are defined as the images of those of $\Delta^n$ with respect to $\varphi$.

The manifold $M^n$ is triangulable, **[7]**. It means that for any $l$, $0 \leq l \leq n$ such a finite set $\Phi^l$ of diffeomorphisms $\varphi : \Delta^l \to M^n$ is defined that

a) $M^n$ is a disjunct union of images $\varphi(Int\Delta^l)$, $\varphi \in \Phi^l$;

b) if $\varphi \in \Phi^l$ then $\varphi \circ \varepsilon_i \in \Phi^{l-1}$ for every $i$ where $\varepsilon_i : \Delta^{k-1} \to \Delta^k$ is the linear mapping transferring the vertexes $v_0,...,v_{k-1}$ of the simplex $\Delta^{k-1}$ in the vertexes $v_0,...,\hat{v}_i,...v_k$ of the simplex $\Delta^k$.

**2°.** Let $\delta_0^n$ be some simplex of the fixed triangulation of the manifold $M^n$. We paint the inner part $Int\delta_0^n$ of the simplex $\delta_0^n$ in white and the boundary $\partial\delta_0^n$ of $\delta_0^n$ in black. There exist coordinates on $Int\delta_0^n$ given by diffeomorphism $\varphi_0$. A subsimplex $\delta_{01}^{n-1} \subset \delta_0^n$ is defined by a black face $\delta_{01}^{n-1} \subset \delta_0^n$ and the center $c_0$ of $\delta_0^n$. We connect $c_0$ with the center $d_0$ of the face $\delta_{01}^{n-1}$ and decompose the subsimplex $\delta_{01}^n$ as a set of intervals which are parallel to the interval $c_0 d_0$. The face $\delta_{01}^{n-1}$ is a face of some simplex $\delta_1^n$ that has not been painted. We draw an interval between $d_0$ and the vertex $v_1$ of the subsimplex $\delta_1^n$ which is opposite to the face $\delta_{01}^{n-1}$ then we decompose $\delta_1^n$ as a set of intervals which are parallel to the interval $d_0 v_1$. The set $\delta_{01}^n \cup \delta_1^n$ is a union of such broken lines every one from which consists of two intervals where the endpoint of the first interval coincides with the beginning of the second interval (in the face $\delta_{01}^{n-1}$) the first interval belongs to $\delta_{01}^n$ and the second interval belongs to $\delta_1^n$. We construct a homeomorphism (extension) $\varphi_{01}^1 : Int\delta_{01}^n \to Int(\delta_{01}^n \cup \delta_1^n)$. Let us consider a point $x \in Int\delta_{01}^n$ and let $x$ belong to a broken line consisting of two intervals the first interval is of a length of $s_1$ and the second interval is of a length of $s_2$ and let $x$ be at a distance of $s$ from the beginning of the first interval. Then we suppose that $\varphi_{01}^1(x)$ belongs to the same broken line at a distance of $\frac{s_1 + s_2}{s_1} \cdot s$ from the beginning of the first interval. It is clear that $\varphi_{01}^1$ is a homeomorphism giving coordinates on $Int(\delta_{01}^n \cup \delta_1^n)$. We paint points of $Int(\delta_{01}^n \cup \delta_1^n)$ white. Assuming the coordinates of points of white initial faces of subsimplex $\delta_{01}^n$ to be fixed we obtain correctly introduced coordinates on $Int(\delta_0^n \cup \delta_1^n)$. The set $\sigma_1 = \delta_0^n \cup \delta_1^n$ is called a *canonical polyhedron*. We paint faces of the boundary $\partial \sigma_1$ black.

We describe the contents of the successive step of the algorithm of extension of coordinate neighborhood. Let us have a canonical polyhedron $\sigma_{k-1}$ with white inner points (they have introduced *white coordinates*) and the black boundary $\partial \sigma_{k-1}$. We look for such an $n$–simplex in $\sigma_{k-1}$, let it be $\delta_0^n$ that has such a black face, let it be $\delta_{01}^{n-1}$ that is simultaneously a face of some $n$–simplex, let it be $\delta_1^n$, inner points of which are not painted. Then we apply the procedure described above to the pair $\delta_0^n$, $\delta_1^n$. As a result we have a polyhedron $\sigma_k$ with one simplex more than $\sigma_{k-1}$ has. Points of $Int\sigma_k$ are painted in white and the boundary $\partial \sigma_k$ is painted in black. The process is finished in the case when all the black faces of the last polyhedron border on the set of white points (the cell) from two sides.

After that all the points of the manifold $M^n$ are painted in black or white, otherwise we would have that $M^n = M_0^n \cup M_1^n$ (the points of $M_0^n$ would be painted and those of $M_1^n$ would be not) with $M_0^n$ and $M_1^n$ being unconnected, which would contradict of connectivity of $M^n$.

Thus, we have proved the following

**Theorem 1.** *Let $M^n$ be a connected, compact, closed, smooth manifold of dimension n. Then $M^n = C^n \cup K^{n-1}$, $C^n \cap K^{n-1} = \emptyset$, where $C^n$ is an n–dimensional cell and $K^{n-1}$ is a union of some finite number of (n–1)–simplexes of the triangulation.*

**3°.** We consider the initial simplex $\delta_0^n$ of the triangulation and its center $c_0$. Drawing intervals between the point $c_0$ and points of all the faces of $\delta_0^n$ we obtain a decomposition of $\delta_0^n$ as a set of the intervals. In **2°** the homeomorphism $\psi : Int\delta_0^n \to C^n$ was constructed and $\psi$ evidently maps every interval above on a piecewise smooth broken line $\gamma$ in $C^n$. We denote $\widetilde{M}^n = M^n \setminus \{c_0\}$. $\widetilde{M}^n$ is a connected and simply connected manifold if $M^n$ is that. Let $I=[0;1]$, we define a homotopy $F: \widetilde{M}^n \times I \to \widetilde{M}^n : (x; t) \mapsto y=F(x;t)$ in the following way

a) $F(z; t)=z$ for every point $z \in K^{n-1}$;
b) if a point $x$ belongs to the broken line $\gamma$ in $C^n$ and the distance between $x$ and its limit point $z \in K^{n-1}$ is $s(x)$ then $y=F(x; t)$ is on the same broken line $\gamma$ at a distance of $(1-t)s(x)$ from the point $z$.

It is clear that $F(x;0)=x$, $F(x;1)=z$ and we have obtained the following

**Theorem 2.** *The spaces $\widetilde{M}^n$ and $K^{n-1}$ are homotopy–equivalent, in particular, the groups of singular homologies $H_k(\widetilde{M}^n)$ and $H_k(K^{n-1})$ are isomorphic for every k.*

**Corollary 2.1.** *The space $K^{n-1}$ is connected and if $M^n$ is simply connected then $K^{n-1}$ is simply connected too.*

**Remark.** *The white coordinates are extended from the simplex $\delta_0^n$ in the simplex $\delta_1^n$ through the face $\delta_{01}^{n-1}$ hence $Int\delta_{01}^{n-1}$ has also the white coordinates. On the other hand there exist two linear structures (intervals, the center etc) on $\delta_{01}^n$ induced from $\delta_0^n$ and $\delta_1^n$ respectively. Further, we set that the linear structure of $\delta_{01}^{n-1}$ is the structure induced from $\delta_0^n$.*

## 2. On the complex $K^2$

For a three–dimensional, connected, compact, closed, smooth manifold $M^3$ we consider a decomposition $M^3 = C^3 \cup K^2$ obtained in theorem 1.

We call simplexes of dimension 3, 2, 1 by tetrahedrons, triangles, edges (intervals) respectively.

**1°. Definition 1.** a) *A triangle from the complex $K^2$ is called a f–triangle (free) if it has at least one free edge i. e. such an edge that it is not an edge of any other triangle from $K^2$.*

b) *A triangle from the complex $K^2$ is called a m–triangle if it has such an edge that is an edge of more than two triangle from $K^2$. A m–triangle can not be a f–triangle.*

c) *A triangle from the complex $K^2$ is called a s–triangle (standard) if every its edge is an edge of only one other triangle from $K^2$.*

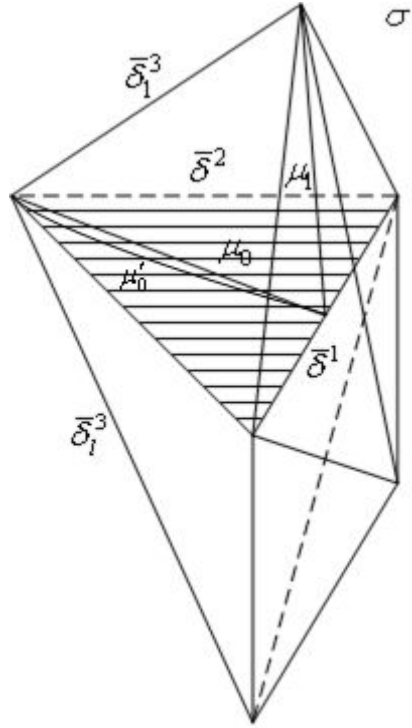

**Fig**. 1.

Let us have a *f*–triangle $\bar{\delta}^2 \in K^2$ with some free edge $\bar{\delta}^1$. We consider such a polyhedron $\sigma$ which $\sigma$ is a set of all the tetrahedrons with $\bar{\delta}^1$ as their edge. Among them we have exactly two tetrahedrons, let they be $\bar{\delta}_1^3$ and $\bar{\delta}_l^3$ with $\bar{\delta}^2$ as their face. We call the *output* of $\bar{\delta}_1^3$ the face $\bar{\delta}_1^2$ with $\bar{\delta}^1$ as its edge. Inner points of the triangle $\bar{\delta}_1^2$ are white because the edge $\bar{\delta}^1$ is free. The face $\bar{\delta}_1^2$ is a face of another tetrahedron $\bar{\delta}_2^3$ that has only one another face $\bar{\delta}_2^2$ with the edge $\bar{\delta}^1$, moreover, all inner points of the triangle $\bar{\delta}_2^2$ are white. The faces $\bar{\delta}_1^2$ and $\bar{\delta}_2^2$ are called respectively the *input* and *output* (*conversions*) of the tetrahedron $\bar{\delta}_2^3$. The face $\bar{\delta}_2^2$ is called the input of some tetrahedron $\bar{\delta}_3^3$ *etc*. Taking a finite number of steps we come to the tetrahedron $\bar{\delta}_l^3$ with an input $\bar{\delta}_{l-1}^2$ with $\bar{\delta}^1$ as its edge and all inner points of the triangle $\bar{\delta}_{l-1}^2$ are white. Thus, we obtain $\sigma = \bigcup_{i=1}^{l} \bar{\delta}_i^3$ (minimal possible meaning is $l=3$). We have to note that all inner points of the faces of conversions $\overline{\delta}_1^2, ..., \overline{\delta}_{l-1}^2$ in the tetrahedrons of the polyhedron $\sigma$ are white.

We call the midpoint of an interval a point that divides this interval into two equal parts by length in metric $g$ (it might not be identical to a center of the interval). As usually, we call a median of a triangle an interval connecting a vertex of the triangle and the midpoint of the opposite side. Such a median is unique in white triangle (see remark from **3°**, 1) there are two such medians $\mu_0$ and $\mu_0'$ in a black triangle. We decompose the tetrahedrons $\overline{\delta}_1^3$ and $\overline{\delta}_l^3$ as sets of intervals every one of them is parallel to the median $\mu_0$ ($\mu_0'$) if the interval belongs to $\overline{\delta}_1^3$ ($\overline{\delta}_l^3$) where the medians $\mu_0$ and $\mu_0'$ of the triangle $\overline{\delta}^2$ (provided with two linear structures) are drawn to the midpoint of the edge $\overline{\delta}$. One of the endpoints of every such an interval belongs to $\overline{\delta}_1^2$ or $\overline{\delta}_{l-1}^2$ and the other one belongs to the boundary of the polyhedron $\sigma$. Further, we decompose every one of the tetrahedrons $\overline{\delta}_2^3, ..., \overline{\delta}_{l-1}^3$ as a set of intervals that are parallel to the edge of the tetrahedron which is not in conversions of the tetrahedron. Such intervals connect the points of the input and the output of the tetrahedron. Thus, the polyhedron $\sigma$ is decomposed as a set of $l$–broken lines. We should mention that $(Int\sigma)\backslash\overline{\delta}^2$ and some part of the boundary of $\sigma$ are painted in white i. e. their points have white coordinates.

**2°. Proposition 3.** *We can redistribute coordinates of white points of the polyhedron $\sigma$ and introduce white coordinates of points from $Int\overline{\delta}^2 \cup \overline{\delta}^1$ (construct the corresponding homemorphism $\varphi_\sigma$) in such way that the following conditions are fulfilled*

a) *all the points of Int $\sigma$ are painted in white i.e. have white coordinates,*

b) *white coordinates of points of boundary faces of the polyhedron $\sigma$ are not changed.*

**Proof.** Let a white point $x \in (Int\sigma)\backslash\overline{\delta}^2$ then it belongs to a broken line $x_0x_1...x_l$, where $x_i \in \overline{\delta}_i^2$, $i=1, ..., l-1$, $x_0$ is the beginning of the corresponding interval from the decomposition of $\overline{\delta}_1^3$, and $x_l$ is the endpoint of the corresponding interval from the decomposition of $\overline{\delta}_l^3$, $x_{i-1}x_i$ is the interval from the decomposition of $\overline{\delta}_i^3$ (see **1°**). We draw the median $\mu_1$ to the edge $\overline{\delta}^1$ of the triangle $\overline{\delta}_1^2$ then we decompose $\overline{\delta}_1^2$ as a set of intervals pavallel to $\mu_1$. We consider broken lines consisting of two intervals of the decompositions of $\overline{\delta}_1^2$ and $\overline{\delta}^2$ in the tetrahedron $\overline{\delta}_1^3$ where the beginning of the other interval coincides with the endpoint of the first one in $\overline{\delta}^1$. Let a point $x_1$ belong to some interval of a length of $s_1$ from the decomposition of $\overline{\delta}_1^2$ at a distance of $s$ from the beginning of the interval, let the other interval of the corresponding broken line is of a length of $s_2$ and lies in $\overline{\delta}^2$. We consider a *mapping by length* $\varphi_\sigma$: $\overline{\delta}_1^2 \to \overline{\delta}_1^2 \cup \overline{\delta}^2$: $x_1 \mapsto y_1 = \varphi_\sigma(x_1)$ where the point $y_1$ belongs to the same broken line at a distance of $\frac{s_1+s_2}{s_1}s$ from the beginning of the first interval. It is clear that $\varphi_\sigma$ is a homeomorphism from $\operatorname{Int}\bar\delta_1^2$ on $\operatorname{Int}(\bar\delta_1^2\cup\bar\delta^2)$ introducing white coordinates on $\operatorname{Int}(\bar\delta_1^2\cup\bar\delta^2)$.

Let us draw an interval between $x_0$ and $y_1$ in $\bar\delta_1^3$ we shall consider the following cases.

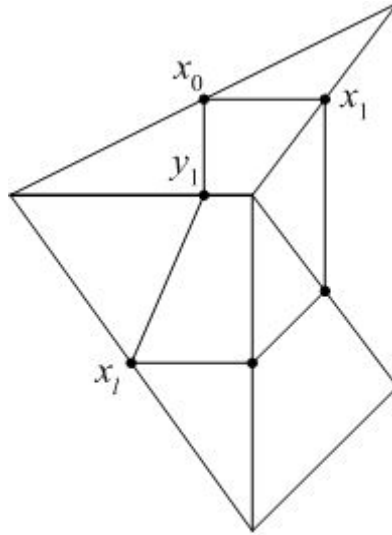

**Fig**. 2.1

1) $y_1\in\bar\delta^2$. We draw the interval $y_1x_l$ in the tetrahedron $\bar\delta_l^3$ and set that $\varphi_\sigma$ maps the broken line $x_0x_1...x_l$ on the broken line $x_0y_1x_l$ by length and $\varphi_\sigma(x_0)=x_0$, $\varphi_\sigma(x_l)=x_l$, $\varphi_\sigma(x_1)=y_1$, therefore the point $y=\varphi_\sigma(x)$, $y\in x_0y_1x_l$, is defined uniquely.

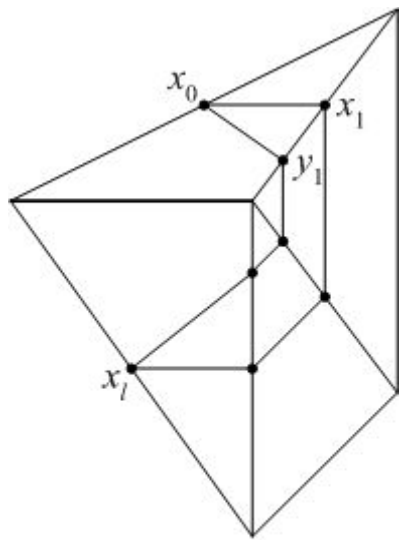

**Fig**. 2.2

2) $y_1 \in \bar{\delta}^2$. We consider the broken line $x_0 y_1 y_2 ... y_{l-1} x_l$ where $y_{i-1} y_i$ is an interval of the decomposition of the tetrahedron $\bar{\delta}_i^3$ constructed in **2°** and we suppose that $\varphi_\sigma$ maps the broken line $x_0 x_1 ... x_l$ on the broken line $x_0 y_1 y_2 ... y_{l-1} x_l$ by length and $\varphi_\sigma(x_0)=x_0$, $\varphi_\sigma(x_l)=x_l$, $\varphi_\sigma(x_i)=y_i$, $i=1, ..., l-1$, therefore the point $y=\varphi_\sigma(x)$, $y \in x_0 y_1 y_2 ... y_{l-1} x_l$ is uniquely defined.

3) For every white point $x$ of the boundary of the polyhedron $\sigma$ we suppose that $\varphi_\sigma(x)=x$.

Thus, one–to–one mapping $\varphi_\sigma: (\sigma \setminus \bar{\delta}^2) \to \sigma$ has been constructed. It is clear that $\varphi_\sigma$ is an homemorphism. Obviously, the conditions a) and b) are fulfilled

**QED.**

The constructed homemorphism $\varphi_\sigma$ makes possible to paint any *f*–triangle white. It is clear that process of painting of *f*–triangles described in the proof of proposition 3 is reduced to retracting every of them to corresponding two edges, therefore, it does not change the simple connectivity of space on every step. Sets of black and white points will be denoted by $K^2$ and $C^3$ respectively on every step.

**3°. Definition 2.** *An edge $\bar{\delta}^1 = x_0 x_l$ is called isolated if it is not an edge of any triangle from $K^2$ and one of the endpoints of the interval $\bar{\delta}^1$ (let it be $x_l$) is free i. e. it is not an endpoint of any edge from $K^2$.*

An isolated edge $\bar{\delta}^1$ can appear as a result of painting in white some neighboring *f*–triangles containing $\bar{\delta}^1$. We consider polyhedrons $\sigma$ and $\bar{\sigma}$ where $\sigma$ is the set of all tetrahedrons with $x_l$ as their vertex and $\bar{\sigma}$ is the set of all tetrahedrons with $\bar{\delta}^1$ as their edge. It is clear that $\bar{\sigma} \subset \sigma$ and all the points of Int$\sigma$ are white with the exception of black points of $\bar{\delta}^1 \setminus \{x_0\}$.

**Proposition 4.** *We can redistribute coordinates of white points of the polyhedron $\sigma$ and introduce white coordinates of points from Int$\bar{\delta}^1 \cup \{x_l\}$ (construct the corresponding homeomorphism) in such a way that the following condition are fulfilled)*

a) *all the points of Int$\sigma$ are painted in white i. e. have white coordinates,*

b) *white coordinates of points of boundary faces of the polyhedron $\sigma$ are not changed.*

**Proof.** We shall divide the proof into two steps.

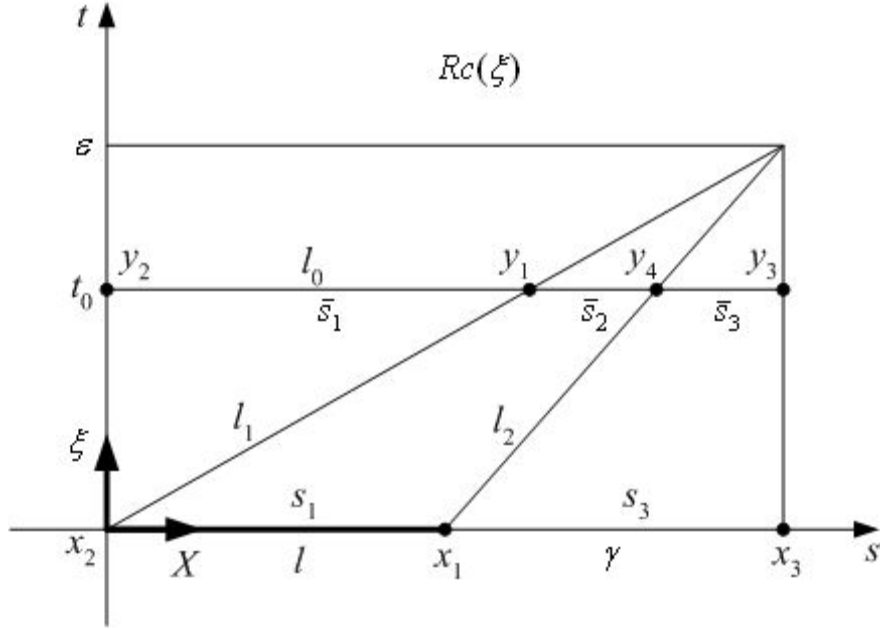

**Fig**. 3.

1) We consider the Riemannian manifold ($M^3$, g). Let a point $x_2$ be the midpoint of the interval $\bar{\delta}^1 = x_0 x$. By the definition of triangulation $x_2 x_1$ is a segment of a differentiable curve $x = x(s)$ of a length of $s_1$ ($x(0) = x_2$, $x(s_1) = x_1$) and the tangent vector $X_x = \dot{x}(s)$ is defined in every point of the curve. The pair ($x_1$, $X_{x_1}$) defines a segment of a geodesic $\gamma$ : $x = x(s)$ where $s \in [0; s_3]$. We can choose a segment $\gamma$ so small that $\gamma$ belongs to $Int\sigma$ and $\gamma$ has no points of self–intersection and points of intersection with $\bar{\delta}^1$. Thus, we can consider the segment of the curve $l$ : $x = x(s)$, $s \in [0; s_1 + s_3]$ consisting of two parts above and belonging to $Int\sigma$. There exists sufficiently small *normal tubular neighborhood* $Tb(l) \subset Int\sigma$, [6]. Further,
$Tb(l) = \bigcup D(x(s), \varepsilon)$, where $D(x(s); \varepsilon) = \{y = \exp_{x(s)}(t\xi) \mid \xi \perp X_{x(s)}, \|\xi\| = 1, 0 \leq t \leq \varepsilon\}$.

For a fixed vector $\xi_{x_2}$, $\|\xi_{x_2}\| = 1$, $\xi_{x_2} \perp X_{x_2}$ by the parallel transport with respect to the Riemannian connection $\nabla$ there exists the unique vector field $\xi$ along the curve $l$ and the rectangle $Rc(\xi) = \{y = \exp_{x(s)}(t\xi) \mid 0 \leq s \leq s_1 + s_3, 0 \leq t \leq \varepsilon\}$. In the rectangle $Rc(\xi)$ we consider the following segments of curves $l_0$: $y = \exp_{x(s)}(t_0 \xi)$, $t = t_0$, $s \in [0; s_1 + s_3]$; $l_1$: $y = \exp_{x(s)}\left(\dfrac{\varepsilon s}{s_1 + s_3}\xi\right)$, $s \in [0; s_1 + s_3]$; $l_2$: $y = \exp_{x(s)}\left(\dfrac{\varepsilon}{s_3}(s - s_1)\xi\right)$, $s \in [s_1; s_1 + s_3]$. Let $y_2 = \exp_{x_2}(t_0 \xi) \in l_0$, $y_3 = \exp_{x_3}(t_0 \xi) \in l_0$, $y_1 = l_0 \cap l_1$, $y_4 = l_0 \cap l_2$ and the corresponding lengthes of segments of the curve $l_0$ are equal to $\bar{s}_1$ for $y_2 y_1$, $\bar{s}_2$ for $y_1 y_4$, $\bar{s}_3$ for $y_4 y_3$. A mapping $\varphi_1$ is defined in the following way
a) $\varphi_{1|\sigma \setminus Tb(l)} = id$; b) if $y \in y_4 y_3 \subset l_0$ and $y$ is at a distance of $s$ from $y_4$ then $\varphi_1(y) \in y_1 y_3 \subset l_0$ and $\varphi_1(y)$ is at a distance of $\dfrac{\bar{s}_1 + \bar{s}_3}{\bar{s}_3} s$ from $y_1$, in particular, $\varphi_1(y_4) = y_1$; c) if $y \in y_2y_4 \subset l_0$ and $y$ is at a distance of $s$ from $y_2$ then $\varphi_1(y) \in y_2y_1 \subset l_0$ and $\varphi_1(y)$ is at a distance of $\dfrac{\bar{s}_1}{\bar{s}_1 + \bar{s}_2}s$ from $y_2$, in particular, $\varphi_1(y_4)=y_1$; d) if for $l_0 \quad t_0=\varepsilon$, then $\varphi_{1|l_0}=id$; e) if $x \in x_1x_3 \subset l$ and $x$ is at a distance of $s$ from $x_1$ then $\varphi_1(x) \in x_2x_3 \subset l$ and $\varphi_1(x)$ is at a distance of $\dfrac{s_1+s_3}{s_3}s$ from $x_2$, in particular, $\varphi_1(x_1)=x_2$. An one–to–one correspondence $\varphi_1 : (Tb(l) \setminus \bar{\delta}^1) \to (Tb(l) \setminus \{x_2\})$ has been constructed and $\varphi_1$ is an homeomorphism introducing white coordinates on $(Tb(l) \setminus \{x_2\})$, in addition, $\varphi_1=id$ on the boundary $\partial Tb(l) \setminus \{x_2\}$. Thus, $Int(x_2x_1) \cup \{x_1\}$ is painted in white.

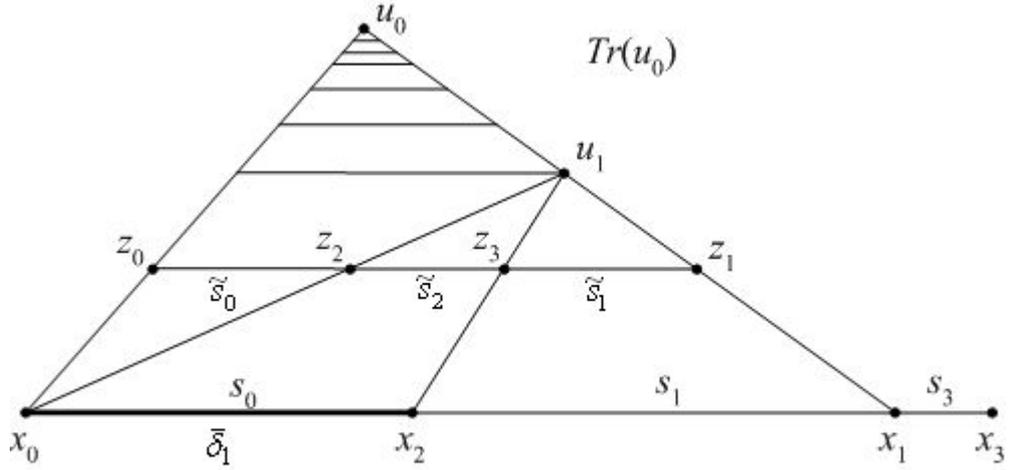

**Fig**. 4.

2)  All the points of $Int\sigma$ are white except for black points of the interval $x_0x_2 \subset \bar{\delta}^1$. Let $\bar{\delta}^3$ be some tetrahedron from $\bar{\sigma}$ then we consider a cross–section $Tr(u_0)$ (triangle) of $\bar{\delta}^3$ passing through the edge $\bar{\delta}^1$ and a point $u_0$ belonging to that edge of $\bar{\delta}^3$ that is crossed with $\bar{\delta}^1$. We draw the median $x_0u_1$ and the interval $u_1x_2$ in the triangle $Tr(u_0)$, where $u_1$ is the midpoint of the interval $u_0x_1$. Let $z_0z_1$ be any interval in $Tr(u_0)$ parallel to $\bar{\delta}^1$ where $z_0 \in u_0x_0$, $z_1 \in u_1x_1$, $z_2$ is the point of intersection of the interval $z_0z_1$ and the median $x_0z_1$, $z_3$ is the point of intersection of the intervals $z_0z_1$ and $u_1x_2$ the lengthes of intervals are respectively equal to $\tilde{s}_0$ for $z_0z_2$, $\tilde{s}_2$ for $z_2z_3$, $\tilde{s}_1$ for $z_3z_1$.

A mapping $\varphi_2$ is defined in the following way

a) $\varphi_2(z_0)=z_0$, $\varphi_2(z_1)=z_1$ and $\varphi_2(z)=z$ for $z \in z_0z_1$ where $z_1 \in u_0u_1$; b) if $z \in z_3z_1$ and $z$ is at a distance of $s$ from $z_3$ then $\varphi_2(z) \in z_2z_1$ and $\varphi_2(z)$ is at a distance of $\dfrac{\tilde{s}_2 + \tilde{s}_1}{\tilde{s}_1}s$ from $z_2$, in particular, $\varphi_2(z_3)=z_2$; c) if $z \in z_0z_3$ and $z$ is at a distance of $s$ from $z_0$ then $\varphi_2(z) \in z_0z_2$ and $\varphi_2(z)$ is at a distance of $\dfrac{\tilde{s}_0}{\tilde{s}_0 + \tilde{s}_2}s$ from $z_0$, in particular, $\varphi_2(z_3)=z_2$; d) if $x \in x_2x_1$ and $x$ is at a distance of $s$ from $x_2$ then $\varphi_2(x) \in x_0x_1$ and $\varphi_2(x)$ is at a distance of $\dfrac{s_0 + s_1}{s_1}s$ from $x_0$, in particular, $\varphi_2(x_2)=x_0$.

It is clear that $\varphi_2=id$ on the boundary $\partial\overline{\sigma}$ and that $\varphi_2$ is a homeomorphism between $\overline{\sigma}\setminus\{x_0x_2\}$ and $\overline{\sigma}\setminus\{x_0\}$, introducing white coordinates on $Int\,\overline{\sigma}$.

Thus, the composition of the homeomorphism $\varphi_1$ and $\varphi_2$ gives us the realization of all the conditions of the proposition.

**QED.**

**4°.** We assume that in the process of painting $f$–triangles white by the proposition 3 all the triangles from $K^2$ are white i.e. that we have a representation $M^3 = C^3 \cup K^1$, $C^3 \cap K^1 = \emptyset$, where $C^3$ is a three–dimensional cell and $K^1$ is a connected union of finite number of black edges of the triangulation. Since the process of painting $f$–triangles white does not influence simple connectedness of a space that is been obtained after every step then $K^1$ is a tree if the complex $K^2$ is simply connected. Painting isolated edges of $K^1$ in white by the proposition 4 as a result we have unique black point $x_0$. Thus, we obtain a representation $M^3 = C^3 \cup B(x_0;\varepsilon)$, where $B(x_0;\varepsilon)$ is an open geodesic ball with the center in $x_0$ and of radius $\varepsilon$. The manifold $M^3$ is homeomorphic to sphere $S^3$ by the following lemma 5.

**Lemma 5 [5].** *If a topological manifold $M^n$ is a union of two n–dimensional cells then $M^n$ is homeomorphic to sphere $S^n$.*

### 3. Proof of the main theorem

The proof has a combinatorial nature and assumes the realization of a number of algorithms. We consider that step by step. The initial complex $K^2$ is assumed to be connected, simply connected and without free triangles.

**1°.** We call a sequence of tetrahedrons (triangles, edges) a *simple chain (s – chain)* if every such a simplex participates in the sequence only one time and if every subsequent tetrahedron (triangle, edge) has a common face (edge, vertex) with the previous one. The number of elements of an s–chain is called the *length* of the s–chain.

Let $\delta_0^2$ and $\delta_1^2$ be two triangles from the complex $K^2$ with $\delta_0^1$ as their common edge. The edge $\delta_0^1$ can also be an edge of some $m$–triangles other than $\delta_0^2$ and $\delta_1^2$.

**Lemma 6.** *We can rebuild the complex $K^2$ in such a way that $\delta_0^2$ and $\delta_1^2$ are standard by the edge $\delta_0^1$ i. e. we have no m–triangles by the edge $\delta_0^1$. A new rebuilt complex $K^2$ is connected and simple connected.*

**Proof.** We consider an s–chain of tetrahedrons with $\delta_0^1$ as their edge the first of which has the *upper part* of $\delta_0^2$ as its face and the last of which has the upper part of $\delta_1^2$ as its face (the notions *upper* and *lower parts* are naturally defined if we consider a small geodesic ball with a center in an inner point of the edge $\delta_0^1$, this ball is divided by the union $\delta_0^2 \cup \delta_1^2$ into two open subsets which are called upper and lower respectively).

In this s–chain we can find a tetrahedron, let it be $\bar{\delta}_1^3$, which is the first from the s–chain to have a *m–triangle* as its face. Thus, we obtain an s–chain $\bar{\delta}_1^3$, ..., $\bar{\delta}_{l_1}^3$ of tetrahedrons (some of them have *m–triangles* as their faces) $\bar{\delta}_{l_1}^3$ has the upper part of $\delta_1^2$ as its face. We cancel the white painting of points of $\bar{\delta}_2^3$, ..., $\bar{\delta}_{l_1}^3$. Then we consider the center of the tetrahedron $\bar{\delta}_1^3$ and realize the of extension of white coordinates through all the faces (including *m–triangles*) of the s–chain $\bar{\delta}_1^3$, ..., $\bar{\delta}_{l_1}^3$ with $\delta_0^1$ as their edge. So, $Int(\bar{\delta}_2^3 \cup ... \cup \bar{\delta}_{l_1}^3)$ and all the inner points of conversions are painted in white, other faces of $\bar{\delta}_2^3$, ..., $\bar{\delta}_{l_1}^3$ (except for conversions) are painted in black.

We shall note that in the process of repainting those faces in black so–called *dead end* might appear. We consider the tree $L$ connecting by intervals all the centers of the tetrahedrons of the triangulation and the centers of all the white faces. If we repaint a face $\bar{\delta}^2$ of one from tetrahedrons $\bar{\delta}^3$ of an s–chain from white in to black the tree $L$ might be divided by $\bar{\delta}^2$ into two parts $\bar{L}$ and $L'$ where $\bar{L}$ contains the center of $\bar{\delta}^3$ and $L'$ defines tetrahedrons forming a dead end. Repainting the dead end in black (its boundary has been painted in this colour) we carry out the extension of white coordinates from $\bar{\delta}^3$ through the face $\bar{\delta}^2$ as it was described in section **1**. As a result all the points of the dead end and inner points of $\bar{\delta}^2$ are painted in white and $L'$ is connected with $\bar{L}$ (the simplest example is a dead end made of one tetrahedron). If we have several dead ends then the procedure above is applied to every one of them. The similar operation can be done from the lower part of $\delta_0^2$ on to the lower part of $\delta_1^2$. As a result of the procedure the edge $\delta_0^1$ is an edge of only two triangles $\delta_0^2$ and $\delta_1^2$ of a new rebuilt connected and simple connected complex $K^2$ i. e. $\delta_0^2$ and $\delta_1^2$ are standard by the edge $\delta_0^1$. It is obvious that the set of all the white points remains a three–dimensional cell.

**QED.**

It is well known (see, for example, [**3**]) that a connected and simple connected complex $K^2$ is homotopy–equivalent to union of two–dimensional spheres. We set a task of rebuilding the complex $K^2$ in order to extract from it a two–dimensional disk (two–dimensional sphere) the points of which will be painted in blue.

We consider the process of the extension of a coordinate neighborhood described in section **1**. We choose any triangle $\delta_0^2$ from $K^2$ as an initial polygon $\sigma_0$. The set $Int\sigma_0$ is painted in blue (inner points of $\delta_0^2$ have coordinates given by the corresponding simplex) and the boundary $\partial\sigma_0$ is painted in red. We choose any red edge $\delta_0^1$ in the polygon $\sigma_0$, then we consider any black triangle $\delta_1^2$ from $K^2$ with $\delta_0^1$ as its edge. Using lemma 6 we rebuild $K^2$ moving all $m$–triangles from $\delta_0^1$ on to the boundary $\partial\sigma_1$ where $\sigma_1=\sigma_0\cup\delta_1^2$. Drawing medians in $\delta_0^2$ and $\delta_1^2$ to $\delta_0^1$ we decompose these triangle as sets of intervals which are parallel to the medians. By the two–dimensional version of theorem 1 we can extend blue coordinates from $Int\sigma_0$ into $Int\sigma_1$ painting $Int\sigma_1$ in blue and the boundary $\partial\sigma_1$ of the obtained polygon in red. We shall note that a chosen linear structure on $\delta_1^2$ is coordinated with the linear structure on $\delta_0^2$ i. e. it has been induced by that of the two tetrahedrons which has similar direction («up» or «down») with the tetrahedron inducing the linear structure on $\delta_0^2$.

Further, we consider the following $k$–th step of extension of the *canonical polygon* $\sigma_{k-1}$ where the inner points of $\sigma_{k-1}$ are blue and the red boundary $\partial\sigma_{k-1}$ is homeomorphic to the circle $S^1$. The blue points can not be limit points of the set of black triangles from $K^2$.

**2°.** Let $Ext\sigma_{k-1}$ be the set of all the black triangles from $K^2$ with at least one red edge. We assume that there exists a triangle $\delta_k^2$ in $Ext\sigma_{k-1}$ with two red edges. Extending blue coordinates through one of them we have a red edge $\delta_k^1$ (*inner edge*) one of the endpoints of which (let it be $x_0$) is not an endpoint of some other red edge. The edge $\delta_k^1$ is the common edge of two blue triangles $\delta_k^2$ and $\bar{\delta}_k^2$ and, generally speaking, there exist $m$–triangles with the edge $\delta_k^1$.

**Lemma 7.** *We can rebuild the complex $K^2$ in such a way that we have no black m–triangles on $\delta_k^1$, as well as no a black triangle with the point $x_0$ as its vertex. A new rebuilt complex $K^2$ is connected and simple connected.*

**Proof.** By analogy with lemma 6 we consider an $s$–chain of tetrahedrons with $\delta_0^1$ as their edge the first of which has the upper part of $\delta_k^2$ as its face and the last of which has the upper part of $\delta_1^2$ as its face. We consider the $s$–subchain of tetrahedrons the first of which has the first (towards passing of the $s$–chain) black $m$–triangle on $\delta_k^1$ as its face, the last of which has the last black $m$–triangle on $\delta_k^1$ as its face. Then we apply the algorithm decribed in the proof of lemma 6. As a result of the process we obtain only one black upper *m*–triangle $\delta_b^2$ and one of its edges (let it be $\delta_b^1$) has $x_0$ as its vertex.

We consider the set $\bar{\sigma}$ of all the blue triangles with $x_0$ as their vertex. If we choose a small geodesic ball with the center in $x_0$ then the ball is devided by $\bar{\sigma}$ into two subsets where the *upper semiball* contains the edge $\delta_b^1$ and the *lower semiball* is the other one. We call a *trace* of a simplex with a vertex or endpoint in the point $x_0$ its intersection with the semiball (or its surface called *hemisphere*).

Further, we consider the set of black triangles (they are not *f*–triangles) with the point $x_0$ as their vertex (exept for $\delta_b^2$) and situated «up» i. e. with traces in the upper semiball. We extract *pyramids* from this set. The trace of the *surface of a pyramid* formed of some black triangles is a *maximal closed oval curve* on the upper hemisphere i. e a curve that devides the upper hemisphere into two parts. Any *exterior* white point of the hemisphere close to the *oval* can be connected with the blue boundary by a white curve and any *interior* white point with respect of the oval cannot. Such ovals can be connected among themselves by segments of black curves (they are traces of black triangles called *partitions*).

Further, we consider one of the pyramids and any *s*–chain of white tetrahedrons with the vertex $x_0$ situated «up», the first of which has the upper part of a blue triangle from $\bar{\sigma}$ as its face and the last of which (the first in the *s*–chain) has a black triangle from the surface of the pyramid as its face. In the set of all possible similar *s*–chains we look for an *s*–chain of the minimal length. In the last tetrahedron $\delta_l^3$ of the *s*–chain we consider the subtetrahedron $\delta_{l_0}^3$ with the center of $\delta_b^3$ as its vertex and the mentioned above black triangle as its face. The latter belongs to tetrahedron $\delta_{l_1}^3$. The inner points of $\delta_{l_1}^3$ are simultaneously inner points of the pyramid. Canceling white painting of those inner points and painting all the faces of $\delta_{l_1}^3$ in black we extend white coordinates from $\delta_{l_0}^3$ into $\delta_{l_1}^3$ through their common face as it was described in section **1** and paint those inner points in white again. A new one more length *s*–chain has been obtained and there exist two black faces with $x_0$ as their vertex in the last tetrahedron $\delta_{l_1}^3$ of this *s*–chain. If we obtain a dead end then we eliminate it by the procedure described in lemma 6. It is clear that $\delta_{l_1}^3$ has no blue edge because otherwise the trace of the edge would be a blue inner point with respect to the oval. Further, we iterate the above algorithm and so on. If at some step of the algorithm *f*– triangles appear then we paint them in white by proposition 3. As it has been noted in the proof of this proposition the painting of boundary points of a polyhedron containing a black *f*– triangle is not changed. It is clear that in the end one black edge has been obtained (the oval has been retracted to a point on the hemisphere). If this edge is an edge of a partition then the partition is a *f*– triangle and we can retract it by proposition 3. Similarly, we act with all the upper pyramids and their partitions. Finally, the triangle $\delta_b^2$ becomes free and we can paint it white by proposition 3. The similar procedure can be realized from below.

It is obvious that the set of all the white points is a three–dimensional cell at every step. It is clear that the last rebuilt complex $K^2$ is connected and simple connected because of a homotopy–equivalence.

**QED.**

**Remark.** *So–called semi–isolated edges can remain at the point $x_0$ (the trace of such an edge on a hemisphere is an isolated black point) connecting the point $x_0$ and black two–dimensional subcomplexes which cannot be connected by a black curve, otherwise we have a contradiction to simple connectivity of $K^2$. If we have isolated edges at the point $x_0$ then we can retract them by proposition 4. Further, a structure consisting of a semi–isolated edge and a black subcomplex joined to it is called a «black flower» growing from the point $x_0$.*

Thus, for every point $x \in Int\sigma_k = Int(\sigma_{k-1} \cup \delta_k^2)$ there exists a sufficiently small geodesic ball $B(x)$ that have no points from black triangles in $B(x)$.

**Lemma 8.** *The red points of the set $Int\sigma_k \cap \delta_k^1$ can be painted in blue i. e. all the points from $Int\sigma_k$ have blue coordinates and $Int\sigma_k$ is a two–dimensional cell.*

**Proof** is a simple modification of proposition 3 on the two–dimensional case (for illustration see figures 2.1, 2.2).

**QED.**

**3°.** At some step of the realization of the above algorithm using lemmas 6, 7, 8 an appearance of a tringle $\delta_k^2$ from $Ext\sigma_{k-1}$ with all the three red vertexes is possible. We extend blue coordinates through a red edge of $\delta_k^2$, having released this edge from $m$–triangles by lemma 6. As a result we obtain such a subset $\bar{\sigma}_k = (\sigma_{k-1} \cup \delta_k^2) \subset K^2 = K_0^2$ that its inner points are blue and the red boundary $\partial\bar{\sigma}_k$ is homeomorphic to a pair of circles (a figure eight curve) with a generic point $x_0$. Let $\gamma_0$ be one of the loops of the figure eight curve homeomorphic to $S^1$ and with the support point $x_0$. $\gamma_0$ is a broken line composed of red edges. It follows from simple connectivity of $K_0^2$ there exists a homotopy

$$\bar{F}_0 : I^2 \to K_0^2 : (s; t) \mapsto \gamma_s(t), \bar{F}_0(0; t)=\gamma_0(t), \bar{F}_0(1; t)=x_0, \bar{F}_0(s; 0)= \bar{F}_0(s; 1)=x_0.$$

Thus, three sides of unit square $I^2$ are mapped by $\bar{F}_0$ into the point $x_0$ and its left side ($s=0$) is mapped onto the broken line $\gamma_0$. It is obvious that the loops of the figure eight curve are homotopic to each other in the set of blue points but they can not be retracted in this set hence $\bar{F}_0(I^2)$ contains points of black triangles which are not triangles of black flowers growing on $\sigma_k$. It is clear because any loop from $K^2$ with a black subloop on a flower growing from some point $x$ and also blue (red) points is homotopic to a loop passing through the point $x$ which has no generic points with this flower.

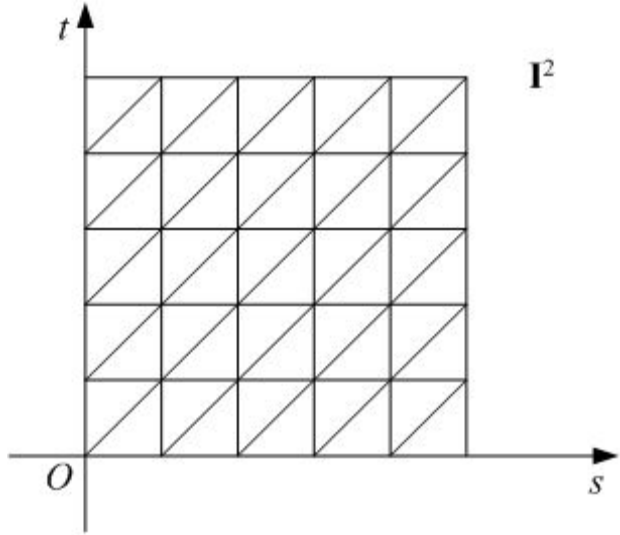

**Fig**. 5.

We mark such points of the left side of the square $I^2$ that their images are vertexes of the broken line $\gamma_0$, then we draw horizontal straight lines through the points. We draw additional horizontal and vertical straight lines and we divide the square $I^2$ into a set of rectangles. Further, every one of the rectangles is divided by its diagonal into two triangles. As a result we obtain so small triangulation $\Sigma_0$ of $I^2$ that the mapping $\overline{F}_0$ has a simplicial approximation $F_0$, [4]. Since $F_0 : I^2 \to K_0^2$ is simplicial, then $F_0$ maps simplexes of the triangulation $\Sigma_0$ onto simplexes of $F_0(I^2)$. Further, we call homotopies only homotopies of the type $F_0$.

**Lemma 9.** *The set of blue points $\overline{\sigma}_k$ bounded by figure eight curve can be reduced to a canonical polyhedron by a retraction of one from its two loops.*

**Proof.** We consider the homotopy $F_0$ and paint simplexes of the square $I^2$ in the colours corresponding to the colours of their images in $F_0(I^2)$. In addition, if two blue triangles from $I^2$ are mapped in the same triangle from $F_0(I^2)$ and if they also have a red simplex (a s–chain of red triangles) as their generic boundary then we repaint this simplex (the s–chain) in blue. If a blue triangle from $I^2$ has a part of the side of $I^2$ (s=0, drawing 5) as its edge or it borders on such an edge through a red simplex (a s–chain of red triangles), then this edge (simplex or s–chain) we also repaint in blue. Thus, any two blue triangles from $I^2$ can be connected by a s–chain of blue triangles and also through a blue edge on the side (s=0) with *exterior* of $I^2$ in the case if there exist blue triangles in $F_0(I^2)$ and out of $F_0(I^2)$. If follows from the fact that the set of blue triangles is indivisible by construction.

Further, it has been noted that $F_0(I^2)$ contains black triangles. We consider a black triangle from $F_0(I^2)$ with the most distant component of the broken line $\gamma_0$ as its edge. Having released this edge from m–triangles by lemma 6 we extend blue coordinates through this edge and we obtain a new broken line $\gamma_1$ in a new complex $K_1^2$. Further, we consider a homotopy $F_1$ mapped the left side of $I^2$ ($s=0$, drawing 5) onto $\gamma_1$ and the other sides of $I^2$ into the point $x_0$.

Out of the black triangles from $F_1(I^2)$ we choose a triangle with the most distant red edge from $x_0$ and repaint this triangle in blue obtaining a new broken line $\gamma_2$ in a new complex $K_2^2$ and so forth. Inner red edges in a complex can be repainted in blue by lemmas 7, 8.

Thus, there exists a sequence of homotopies $F_0$, $F_1$, $F_2$, ... . At every step a set of blue triangles is increased by one triangle. Since a number of such triangles is restricted by a number of all the 2–simplexes of the triangulation then at some step we have a situation when a red broken line $\gamma$ is the boundary of a single black triangle $\delta^2$ with the vertex $x_0$. Repainting inner points of $\delta^2$ in blue and eliminating red edges by lemma 7, 8 we obtain a canonical polygon.

We shall mention that if there are no black triangles in the square then we have no a closed broken line consisting of red simplexes (respectively red loops in $F(I^2)$) in $I^2$ since otherwise we get a contradiction to indivisibility of blue points of the complex (there are not $m$– triangles on blue edges).

**QED.**

**4°.** Using lemmas 6, 7, 8, 9 we continue the extension of a canonical polygon. In addition, the following situation (*a situation of disk*) when at some step a canonical polygon $\sigma_k$ appears is possible, with a furthest extension of $\sigma_k$ being impossible *i. e. Ext* $\sigma_k$ is empty. Black pyramids might exist in the vertexes of a red broken line that is the boundary of $\sigma_k$. By analogy with the process considered in the proof of lemma 7 we transform all those pyramids and partitions (we need not introduce upper and lower pyramids because a set of blue points in a ball forms a *sector of the disk* and the trace of the set on the sphere is *an arc of a circle*) probably, in addition, obtaining black flowers. In this case all the blue triangles with red edges are free. Moving in reverse order to the process of the extension of $\sigma_k$ we consistently contort all the blue triangles by proposition 3. As a result we obtain a new complex of black triangles $\overline{K}^2$ consisting of the black triangles from the initial complex $K^2$ which have no generic points with the polygon $\sigma_k$. Indeed, the process of rebuilding of $K^2$ by lemmas 6, 7 touches upon only those black triangles which have generic points with $\sigma_k$ (in any small ball). Thus, the number of black triangles of the complex $\overline{K}^2$ is certainly less than the number of those in $K^2$. It is obvious that $\overline{K}^2$ is connected and simple connected because of a homotopy–equivalence. By setting $\overline{K}^2$ as an initial complex $K^2$ we begin the construction of a canonical blue polygon again (see **1°**, **2°**, **3°**) and so forth. If the described above situation of a disk is always repeated then et some step of our algorithm the set of black triangles must be exhausted *i. e.* we come to **4°**, **2**.

**Remark.** *The obtained complex can be imagined as a «tree with flowers» growing in the endpoints of the branches of the tree. An iteration of the algorithm can be interpreted as a sequential transformation of those flowers into branches to get a tree in the end.*

**5°.** Using lemmas 6, 7, 8, 9 we continue the extension of a canonical polygon. In the end we have only two variants: the situation of a disk (see **4°**) or a so–called *situation of sphere* considered below. We assume that at some step of our algorithm we have a canonical polygon and its red boundary is the boundary of a black triangle. Extending coordinates inside of this triangle and painting its interior points in blue we apply lemmas 6, 7, 8. As a result we have unique red point $x_0$. It is obvious that the set $M^2$ of all the blue points and the red point $x_0$ is homeomorphic to sphere $S^2$ and $M^2$ is an embedded sectionally smooth submanifold in $M^3$ (in general, smoothness is violated in the points of $K^1$ where $K^1$ is the set of all the edges of triangles forming $M^2$). Let $Tb(K^1) = \bigcup_{x \in K^1} \overline{B}(x; \varepsilon)$ where $\overline{B}(x; \varepsilon)$ is the closed geodesic ball with the center in the point $x$ and of the radius $\varepsilon > 0$. We cancel the white painting of the points of $Tb(K^1)$. It is clear that we can choose such a small $\varepsilon > 0$ that a set of white points forms a 3–dimensional cell as before, further, such an $\varepsilon$ will be considered. Taking into account the smoothness of edges of $K^1$ (any such an edge $l$ has a normal tubular neighborhood $Tb(l)$ considered in proposition 4) it is clear that the surface $M^2$ can be approximated by a smooth surface $\overline{M}^2$ in the following sense

$$M^2 \setminus Tb(K^1) = \overline{M}^2 \setminus Tb(K^1), \quad \overline{M}^2 \subset (M^2 \cup Tb(K^1)).$$

For example, such a procedure can be realized with the help of partition of unity along every edge and in vertexes from $K^1$. Corresponding continuous functions are approximated by smooth functions in the neighborhood $Tb(l)$ of every edge $l$ i. e. *interfaces* of triangles are smoothed. The surface $\overline{M}^2$ is a smooth embedded submanifold in $M^3$.

**Theorem [6].** *Let $M$ be a smooth simple connected manifold and $N$ be its smooth compact connected submanifold of codimension 1. If $\partial M = \partial N = \varnothing$ then $N$ separates $M$.*

So, since $M^3$ is simple connected therefore the submanifold $\overline{M}^2$ must separate $M^3$. But the set of white points is not separated. The obtained contradiction shows that the situation of sphere is impossible.

The main theorem is completely proved.